\documentclass[aps showpacs, showkeys, citeautoscript]{revtex4}
\usepackage{amsfonts}
\usepackage{amssymb}
\usepackage{amstext}
\usepackage{amsmath}
\usepackage{latexsym}
\usepackage{color}
\usepackage{graphicx,color, wasysym}

\newtheorem{defi}{Definition}[section]
\newtheorem{theorem}[defi]{Theorem}

\newcommand{\betheo}{\begin{theorem}$\!\!${\bf \,\,\,}}
\newcommand{\entheo}{\end{theorem}}

  \def\ulamek#1#2{\mbox{\normalfont$\frac{#1}{#2}$}}

\catcode `\@=11 \@addtoreset{equation}{section}

\catcode `\@=12

\begin{document}

\title[On the Laplace transform of the Fr\'{e}chet distribution]
{On the Laplace transform of the Fr\'{e}chet distribution}

\author{K.~A.~Penson}
\email{penson@lptl.jussieu.fr}

\affiliation{Sorbonne Universit\'{e}s, Universit\'e Pierre et Marie Curie (Paris 06), CNRS UMR 7600, \\
Laboratoire de Physique Th\'eorique de la Mati\`{e}re Condens\'{e}e (LPTMC), \\
Tour 13 - 5i\`{e}me \'et., B.C. 121, 4 pl. Jussieu, F 75252 Paris Cedex 05, France,\vspace{2mm}}

\author{K.~G\'{o}rska}
\email{katarzyna.gorska@ifj.edu.pl}

\affiliation{H. Niewodnicza\'{n}ski Institute of Nuclear Physics, Polish Academy of Sciences, Division of Theoretical Physics, ul. Eliasza-Radzikowskiego 152, PL 31-342 Krak\'{o}w, Poland}



\begin{abstract}
We calculate exactly the Laplace transform of the Fr\'{e}chet distribution in the form $\gamma x^{-(1+\gamma)} \exp(-x^{-\gamma})$, $\gamma > 0$, $0 \leq x < \infty$, for arbitrary rational values of the shape parameter $\gamma$, i.e. for $\gamma = l/k$ with $l, k = 1,2, \ldots$. The method employs the inverse Mellin transform. The closed form expressions are obtained in terms of Meijer G functions and their graphical illustrations are provided. A rescaled Fr\'{e}chet distribution serves as a kernel of Fr\'{e}chet integral transform. It turns out that the Fr\'{e}chet transform of one-sided L\'{e}vy law reproduces the Fr\'{e}chet distribution. 

\end{abstract}

\maketitle


\section{Introduction}

We shall be interested in this work in the so-called normalized Fr\'{e}chet distribution in the form 
\begin{equation}\label{eq1}
Fr(\gamma, x) = \gamma x^{-(1+\gamma)} \exp(-x^{-\gamma}), \quad 0 \leq x < \infty,
\end{equation}
for the real shape parameter $\gamma > 0$. It has been introduced by Fr\'{e}chet in 1927 [\onlinecite{MFrechet27}] and has been recognized as one of three types of probability distribution functions (pdf) characterizing extreme value phenomena [\onlinecite{SKotz00}]. The field of applications of Fr\'{e}chet distribution is immense and the reader should consult [\onlinecite{SKotz00}] for references on its role in probability and statistics. Recently the Fr\'{e}chet distribution has been evoked and studied in [\onlinecite{TSimon14}] because of its apparent similarity to one-sided L\'{e}vy stable pdf's which are playing an increasingly important role in non-equilibrium Statistical Mechanics [\onlinecite{KAPenson10, KGorska12, KGorska12a, EBarkai}]. The common features shared between $Fr(\gamma, x)$ and its L\'{e}vy stable counterparts include the essential singularity for $x\to 0$, unimodality and heavy-tail algebraic decay for $x~\to~\infty$. We shall return to this comparison below. 

Since $Fr(\gamma, x)$ is an elementary function, expectation values for many functions on $(0, \infty)$ can be obtained exactly. However, this does not apply to the Laplace transform $\mathbb{E}(e^{-px})$, $\Re(p)>0$ which is explicitly known only in one case $\gamma = 1$ [\onlinecite{TSimon14}]. The knowledge of the Laplace transform is an important issue as it is linked to the alternating moment-generating function.

The purpose of this note is to present exact calculation of the Laplace transform of $Fr(\gamma, x)$ for arbitrary rational values of the shape parameter $\gamma$, i.e. for $\gamma = l/k$ with $l, k = 1,2, \ldots$. The paper is organized as follows: In section II we present the definitions and preliminaries on integral transforms that will serve to obtain the desired results. In section III we give the derivation of the main result and of some relations satisfied by the Laplace transform. We furnish also graphical representations for selected cases. In section IV we introduce the Fr\'{e}chet transform and discuss some of its properties. We conclude the paper in section V where we also give more elements of comparison with L\'{e}vy stable pdf's. 

\section{Definitions and preliminaries}

In the following we give some definitions and informations about the Mellin transform of a function $f(x)$ defined for $x\geq 0$. The Mellin transform is defined for complex $s$ as [\onlinecite{INSneddon72}]
\begin{equation}\label{27.10.13-1}
\mathcal{M}[f(x); s] = f^{\star}(s) = \int_{0}^{\infty} x^{s-1} f(x) dx,
\end{equation}
along with its inverse
\begin{equation}\label{27.10.14-2}
\mathcal{M}^{-1}[f^{\star}(s); x] = f(x) = \frac{1}{2\pi i} \int_{c-i\infty}^{c+i\infty} x^{-s} f^{\star}(s) ds.
\end{equation}
For the role of constant $c$ consult [\onlinecite{INSneddon72}]. For fixed $a>0$, $h\neq 0$, the Mellin transform satisfies the following scaling property:
\begin{equation}\label{27.10.13-3}
\mathcal{M}[x^{b}f(ax^{h}); s] = \frac{1}{|h|} a^{-\frac{s+b}{h}} f^{\star}(\ulamek{s+b}{h}).
\end{equation}
The knowledge of the Mellin transform permits is many cases to obtain the Laplace transform $\mathcal{L}[f(x); p]~=~\int_{0}^{\infty} e^{-px}f(x) dx$, $\Re(p)>0$. This is due to a powerful relation between the Laplace and inverse Mellin transforms:
\begin{equation}\label{27.10.13-4}
\mathcal{M}^{-1}[f^{\star}(1-s)\Gamma(s); x] = \mathcal{L}[f(t); x],
\end{equation}
which belongs to a set of convolution-type theorems discussed in [\onlinecite{LDebnath07}], and is mentioned in [\onlinecite{INSneddon72}], see the  exercise 4-2 (a) on p. 288. The proof of \eqref{27.10.13-4} is immediate as it is sufficient to calculate $\mathcal{M}[\mathcal{L}[f(t); x]; s]$:
\begin{align}\label{27.10.13-5}
\mathcal{M}[\mathcal{L}[f(t); x]; s] &= \mathcal{M}\left[\int_{0}^{\infty} e^{-xt}f(t) dt; s\right] = \int_{0}^{\infty} x^{s-1}\left(\int_{0}^{\infty} e^{-xt}f(t)dt\right)ds \nonumber \\
& = \int_{0}^{\infty} f(t) \Gamma(s) t^{-s} dt = \Gamma(s) \int_{0}^{\infty} t^{-s} f(t) dt = \Gamma(s) f^{\star}(1-s).
\end{align}
The relation \eqref{27.10.13-5} has been previously used in connection with Stieltjes moment problems related to lognormal distributions [\onlinecite{KAPenson99}]. Among known Mellin transforms a very special role is played by those entirely expressible through products and rations of Euler's gamma functions. The Meijer G function is defined as an inverse Mellin transform [\onlinecite{APPrudnikov-v3}]:
\begin{align}\label{27.10.13-6}
G^{m, n}_{p, q}\left(z\Big\vert{\alpha_{1} \ldots \alpha_{p} \atop \beta_{1} \ldots \beta_{q}}\right) &= \mathcal{M}^{-1}\left[\frac{\prod_{j=1}^{m}\Gamma(\beta_{j}+s)\, \prod_{j=1}^{n}\Gamma(1-\alpha_{j}-s)}{\prod_{j=m+1}^{q}\Gamma(1-\beta_{j}-s)\, \prod_{j=n+1}^{p}\Gamma(\alpha_{j}+s)}\right] \\ \label{27.10.13-7}
& = {MeijerG}([[\alpha_{1}, \ldots, \alpha_{n}], [\alpha_{n+1}, \ldots, \alpha_{p}]], [[\beta_{1}, \ldots, \beta_{m}], [\beta_{m+1}, \ldots, \beta_{q}]], z),
\end{align}
where in Eq. \eqref{27.10.13-6} empty products are taken to be equal to 1. In Eqs. \eqref{27.10.13-6} and \eqref{27.10.13-7} the parameters are subject of conditions:
\begin{align}\label{27.10.13-8}
&z\neq 0, \quad 0 \leq m \leq q, \quad 0 \leq n \leq p, \nonumber\\
&\alpha_{j} \in \mathbb{C}, \quad j=1,\ldots, p; \quad \beta_{j}\in\mathbb{C}, \quad j=1,\ldots, q.
\end{align}
For a full description of integration contours in Eq. \eqref{27.10.13-6}, general properties and special cases of the Meijer G functions, see Ref. [\onlinecite{APPrudnikov-v3}]. In Eq. \eqref{27.10.13-7} we present a transparent notation, which we will use in parallel henceforth, inspired by computer algebra systems [\onlinecite{CAS}]. We quote for reference the Gauss-Legendre multiplication formula for gamma functions encountered in this work:
\begin{align}\label{27.10.13-9}
&\Gamma(nz) = (2\pi)^{\frac{1-n}{2}} n^{nz-\frac{1}{2}} \prod_{j=0}^{n-1}\Gamma\left(z+\frac{j}{n}\right), \nonumber\\
&z\neq 0, -1, -2, \ldots, \quad n=1, 2, \ldots.
\end{align}
A frequently occurring special list of $k$ elements in form of
\begin{equation}\label{27.10.13-10}
\frac{a}{k}, \frac{a+1}{k}, \ldots, \frac{a + k-1}{k}
\end{equation}
is denoted by $\Delta(k, a)$, $k\neq 0$, see [\onlinecite{APPrudnikov-v3}].

\section{Derivation of the main result}

It turns out that the Laplace transform of $Fr(\gamma, x)$ can be obtained for $\gamma=l/k$ with $l, k = 1, 2,\ldots$ and the appropriate pdf will be denoted now by $Fr(l, k, x)$. As a first step we calculate the Mellin transform $\mathcal{M}[Fr(l, k, x); s]$, using the property \eqref{27.10.13-3}:
\begin{equation}\label{eq2}
\mathcal{M}[Fr(l, k, x); s] = \Gamma\left[k(\ulamek{1-s}{l} + \ulamek{1}{k})\right],
\end{equation}
which we rewrite as $\Gamma(kz)$ with $z = \ulamek{1-s}{l} + \ulamek{1}{k}$, $z \neq 0, -1, -2, \ldots$. To this last expression we apply the Gauss-Legendre formula \eqref{27.10.13-9} and \eqref{eq2} transforms into 
\begin{equation}\label{eq3}
\Gamma(kz) = \frac{k^{k/l + 1/2}}{(2\pi)^{(k-1)/2}} (k^{k})^{-s/l} \prod_{p=1}^{k} \Gamma\left(\frac{1}{l} + \frac{p}{k} - \frac{s}{l}\right) = \mathcal{M}[Fr(l, k, x); s].
\end{equation}
As a next step we calculate $\mathcal{M}[Fr(l, k, x); 1-s]$ and in order to apply the formula \eqref{27.10.13-4} we apply once again \eqref{27.10.13-9} to $\Gamma(s)$ and obtain:
\begin{equation}\label{eq4}
\Gamma(s) = (2\pi)^{(1-l)/2} l^{-1/2} (l^{l})^{s/l} \prod_{j=1}^{l}\Gamma\left(\frac{s}{l} + \frac{j-1}{l}\right), \quad l=1, 2, \ldots.
\end{equation}
We are in the position now to rewrite
\begin{align}\label{eq5}
\begin{split}
\mathcal{M}[Fr(l, k, x); 1-s] \Gamma(s) &= \frac{k^{1/2} l^{-1/2}}{(2\pi)^{(k+l)/2-1}} (k^{k} l^{l})^{s/l} \left[\prod_{r=1}^{k}\Gamma\left(\frac{s}{l} + \frac{r}{k}\right)\right]\, \left[\prod_{j=1}^{l} \Gamma\left(\frac{s}{l} + \frac{j-1}{l}\right)\right].
\end{split}
\end{align}
Observe that in Eq. \eqref{eq5} the variable $s$ appears only in the combination of $s/l$. Consequently, we invert Eq. \eqref{eq5} with a special case of Eq. \eqref{27.10.13-3}:
\begin{equation}\label{eq6}
\mathcal{M}^{-1}\left[\frac{1}{a^{s/l}} f^{\star}\left(\frac{s}{l}\right); x\right] = l f(a x^{l}), \quad a = (k^{k} l^{l})^{-1}.
\end{equation}

This last expression permits one to use Eqs. \eqref{27.10.13-6} and \eqref{27.10.13-7} and we obtain the final result for the Laplace transform of $Fr(l, k, x)$:
\begin{align}\label{eq7}
\begin{split}
\mathcal{L}[Fr(l, k, x); p] & = \mathcal{M}^{-1}\left[\mathcal{M}[Fr(l, k, x); 1-s] \Gamma(s); p\right] \\
& = \frac{\sqrt{kl}}{(2\pi)^{(k+l)/2 - 1}} G^{k+l, 0}_{0, k+l}\left(\frac{p^{l}}{k^{k} l^{l}} \Big\vert {- \atop \Delta(k, 1), \Delta(l, 0)}\right) \\
& = \frac{\sqrt{kl}}{(2\pi)^{(k+l)/2-1}} MeijerG\left([[\,\,\,], [\,\,\,]], [[\{\ulamek{r}{k}\}_{r=1}^{k}, \{\ulamek{j-1}{l}\}_{j=1}^{l}], [\,\,\,]], \ulamek{p^{l}}{k^{k} l^{l}}\right),
\end{split}
\end{align}
where $\Re(p)>0$. We note that $\{\ulamek{r}{k}\}_{r=1}^{k} = \ulamek{1}{k}, \ulamek{2}{k}, \ulamek{3}{k}, \ldots, 1 = \Delta(k, 1)$ and similarly $\{\ulamek{j-1}{l}\}_{j=1}^{l}=\Delta(l, 0)$, compare the definition Eq. \eqref{27.10.13-10}. The final result in the computer algebra systems notation takes a particularly transparent form:
\begin{equation}\label{eq8}
\mathcal{L}[Fr(l, k, x); p] = \frac{\sqrt{kl}}{(2\pi)^{(k+l)/2 -1}} MeijerG\left([[\,\,\,], [\,\,\,]], [[\Delta(k, 1), \Delta(l, 0)], [\,\,\,]], \ulamek{p^{l}}{k^{k} l^{l}}\right), \quad \Re(p) > 0.
\end{equation}
There is a somewhat hidden symmetry in Eqs. \eqref{eq7} and \eqref{eq8} which becomes apparent by rewriting the list of $k+l$ terms in the third bracket of \eqref{eq8} as
\begin{align}\label{eq9}
\begin{split}
\Delta(k, 1), \Delta(l, 0) & = [\ulamek{1}{k}, \ulamek{2}{k}, \ldots, \ulamek{k-1}{k}], 0, 1, [\ulamek{1}{l}, \ulamek{2}{l}, \ldots, \ulamek{l-1}{l}]\\
& = [\ulamek{1}{l}, \ulamek{2}{l}, \ldots, \ulamek{l-1}{l}], 0, 1, [\ulamek{1}{k}, \ulamek{2}{k}, \ldots, \ulamek{k-1}{k}] \\
& = \Delta(l, 1), \Delta(k, 0).
\end{split}
\end{align}
This last observation leads to a neat transformation law of $\mathcal{L}[Fr(l, k, x); p]$ under the transmutation $l\leftrightarrow k$: 
\begin{equation}\label{eq10}
\mathcal{L}[Fr(l, k, x); p] = \mathcal{L}[Fr(k, l, x); p^{l/k}], \quad \Re(p)>0,
\end{equation}
which links the case $l/k<1$ with that of $l/k>1$. For $l=k=1$ the appropriate Meijer G function can be related to the modified Bessel function $K_{1}(z)$, as then
\begin{align}\label{eq11}
\begin{split}
\mathcal{L}[Fr(1, 1, x); p] &= G^{2, 0}_{0, 2}\left(p\Big\vert {- \atop 1, 0}\right)\\
& = MeijerG([[\,\,\,], [\,\,\,]], [[0, 1], [\,\,\,]], p) \\
& = 2\sqrt{p} K_{1}(2\sqrt{p}), \quad \Re(p) > 0
\end{split}
\end{align}
which is the consequence of formula 8.4.23.1 of {\onlinecite{APPrudnikov-v3}}. The above Laplace transform is the only case that can be represented by a standard special function. Below we list some transforms for several low-lying values of $l, k \geq 1$, with the condition $\Re(p)>0$ applying to all of them:
\begin{align}\label{eq12}
\mathcal{L}[Fr(1, 2, x); p] &= \pi^{-1/2} G^{3, 0}_{0, 3} \left(\frac{p}{4}\Big\vert {- \atop 0, \ulamek{1}{2}, 1}\right) = \pi^{-1/2} MeijerG([[\,\,\,], [\,\,\,]], [[0, \ulamek{1}{2}, 1], [\,\,\,]], \ulamek{p}{4}),\\
\label{eq13}
\mathcal{L}[Fr(2, 1, x); p] &= \pi^{-1/2} G^{3, 0}_{0, 3} \left(\frac{p^{2}}{4}\Big\vert {- \atop 0, \ulamek{1}{2}, 1}\right) = \pi^{-1/2} MeijerG([[\,\,\,], [\,\,\,]], [[0, \ulamek{1}{2}, 1], [\,\,\,]], \ulamek{p^{2}}{4}),
\end{align}
\begin{align}
\mathcal{L}[Fr(1, 3, x); p] &= \frac{\sqrt{3}}{2\pi} G^{4, 0}_{0, 4} \left(\frac{p}{27}\Big\vert {- \atop 0, \ulamek{1}{3}, \ulamek{2}{3}, 1}\right) \nonumber\\
&=  \frac{\sqrt{3}}{2\pi} MeijerG([[\,\,\,], [\,\,\,]], [[0, \ulamek{1}{3}, \ulamek{2}{3}, 1], [\,\,\,]], \ulamek{p}{27}), \label{eq14} \\[0.6\baselineskip]
\mathcal{L}[Fr(2, 3, x); p] &= \frac{\sqrt{12}}{4\pi^{3/2}} G^{5, 0}_{0, 5} \left(\frac{p^{2}}{108}\Big\vert {- \atop 0, \ulamek{1}{3}, \ulamek{1}{2}, \ulamek{2}{3}, 1}\right) \nonumber\\ 
&=  \frac{\sqrt{12}}{4\pi^{3/2}} MeijerG([[\,\,\,], [\,\,\,]], [[0, \ulamek{1}{3}, \ulamek{1}{2}, \ulamek{2}{3}, 1], [\,\,\,]], \ulamek{p^{2}}{108}), \label{eq15} \\[0.6\baselineskip]
\mathcal{L}[Fr(3, 2, x); p] &= \frac{\sqrt{12}}{4\pi^{3/2}} G^{5, 0}_{0, 5} \left(\frac{p^{3}}{108}\Big\vert {- \atop 0, \ulamek{1}{3}, \ulamek{1}{2}, \ulamek{2}{3}, 1}\right) \nonumber\\ 
& = \frac{\sqrt{12}}{4\pi^{3/2}} MeijerG([[\,\,\,], [\,\,\,]], [[0, \ulamek{1}{3}, \ulamek{1}{2}, \ulamek{2}{3}, 1], [\,\,\,]], \ulamek{p^{3}}{108}).\label{eq16}
\end{align}
The reason why the expressions \eqref{eq12}, \eqref{eq13}, \eqref{eq14}, \eqref{eq15} and \eqref{eq16}, as well as the corresponding ones for $l, k~>~3$ (not cited here), cannot be represented by more conventional special functions (as for example the generalized hypergeometric functions ${_{p} F_{q}}$), is that from the list $[\Delta(k, 1), \Delta(l, 0)]$ one can form the differences of indices (members of this lists) that are equal to zero or an integer (positive or negative). This excludes the application of decoupling formulas 8.2.2.3 of [\onlinecite{APPrudnikov-v3}]. Therefore the forms \eqref{eq12}-\eqref{eq16} are the final ones. (The reader might be tempted to use the tabulated integral in the formula 2.3.2.14 on p.~322 of [\onlinecite{APPrudnikov-v1}] which at first sight could correspond to the required Laplace transform. However, this formula is not appropriate as it applies only for $\alpha \neq nr$, $n=0, \pm 1, \pm 2, \ldots$, as can be seen from the arguments of gamma functions.)

The Meijer G functions are implemented in computer algebra systems [\onlinecite{CAS}] and should be considered as full-fledged special functions allowing differentiation, integration etc, as well as plotting.

\begin{figure}[!h]
\begin{center}
\includegraphics[scale=0.56]{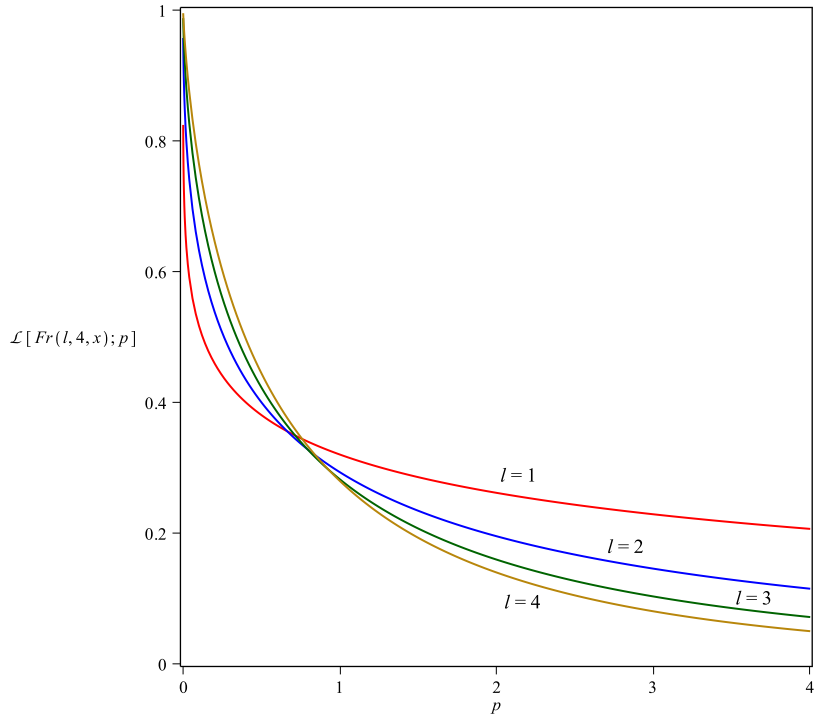}
\caption{\label{fig1} (Color Online) Plot of Eq. \eqref{eq8} for $k=4$, $l=1$ (red line), $l=2$ (blue line), $l=3$ (green line), and $l=4$ (gold line).}
\end{center}
\end{figure}

\begin{figure}[!h]
\begin{center}
\includegraphics[scale=0.56]{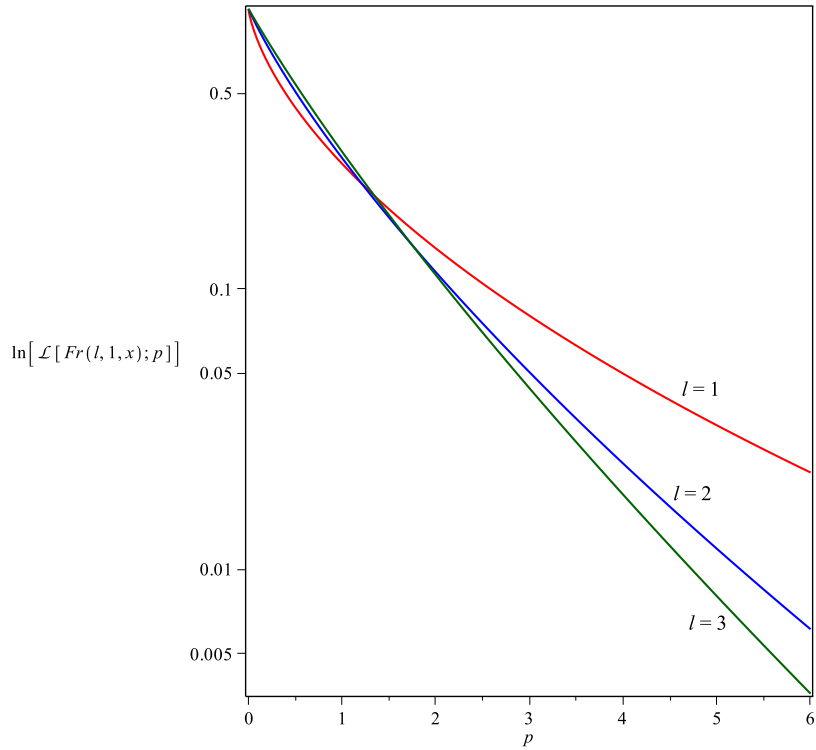}
\caption{\label{fig2} (Color Online) Semi-logarithmic plot of Eq. \eqref{eq8} for $k=1$, $l=1$ (red line), $l=2$ (blue line), and $l=3$ (green line).}
\end{center}
\end{figure}

We have illustrated our results of Eq. \eqref{eq8} by representing graphically the Laplace transforms for selected values of $l$ and $k$. This is done in Figs. \ref{fig1} and \ref{fig2}. Observe that for all $l$ and $k$, $\lim\limits_{p\to 0}\mathcal{L}[Fr(l, k, x); p] = 1$.

\section{Fr\'{e}chet transform}

The knowledge of the Laplace transform of $Fr(\gamma, x)$ permits one to consider a new type of integral transform with the function related to $Fr(\gamma, x)$ as a kernel. In this section we introduce the Fr\'{e}chet integral transform. In an analogy to the L\'{e}vy kernel [\onlinecite{KGorska12a}] the Fr\'{e}chet kernel $\sigma_{\gamma}(x, t)$ is equal to a suitably rescaled two-variable Fr\'{e}chet pdf, namely:
\begin{equation}\label{30.10.13-1}
\sigma_{\gamma}(x, t) = \frac{1}{t^{1/\gamma}} Fr\left(\gamma, \frac{x}{t^{1/\gamma}}\right) = \gamma t x^{-(1+\gamma)} e^{-t x^{-\gamma}}, \quad \gamma, t, x, > 0,
\end{equation}
which is a positive function. For a function $f(x)$ given for $x > 0$, we define its Fr\'{e}chet transform as
\begin{equation}\label{30.10.13-2}
\bar{f}(\gamma, x) = \int_{0}^{\infty} \sigma_{\gamma}(x, t) f(t) dt.
\end{equation}
Since the kernel $\sigma_{\gamma}(x, t)$ contains the exponential function, the Fr\'{e}chet transform is related to the Laplace transform. Indeed, in general we can write
\begin{align}\label{17Feb2014-1}
\bar{f}(\gamma, x) &= \int_{0}^{\infty} \sigma_{\gamma}(x, t) f(t) dt  = \int_{0}^{\infty} \gamma t x^{-(1+\gamma)} e^{-t x^{-\gamma}} f(t) dt \\
& = -\gamma x^{-(1+\gamma)} \int_{0}^{\infty} \left[\frac{d}{d u} e^{-t u}\right]_{u=x^{-\gamma}}\!\! f(t) dt \\
& = -\gamma x^{-(1+\gamma)} \left[\frac{d}{d u} \int_{0}^{\infty} e^{-t u} f(t) dt  \right]_{u=x^{-\gamma}} \\
& = -\gamma x^{-(1+\gamma)} \left\{\frac{d}{d u}  \mathcal{L}[f(t), u]\right\}_{u=x^{-\gamma}}.
\end{align}

For many functions their Fr\'{e}chet transforms can be obtained readily. Among them the one-sided L\'{e}vy stable pdf's $g_{\alpha}(x)$ play a special role. The one-sided L\'{e}vy stable pdf for $0 < \alpha < 1$ by definition satisfies $\mathcal{L}[g_{\alpha}(x); p] = \exp(-p^{\alpha})$, see [\onlinecite{KAPenson10, KGorska12a, EBarkai, Zolotariev-b, Uchaikin-b}]. The Fr\'{e}chet transform of one-sided L\'{e}vy stable distribution $g_{\alpha}(t)$ is calculated below, with $0 < \alpha < 1$:
\begin{align}\label{30.10.13-3}
\bar{g}_{\alpha}(\gamma, x) & = \int_{0}^{\infty} \sigma_{\gamma}(x, t) g_{\alpha}(t) dt = \int_{0}^{\infty} \gamma t x^{-(1+\gamma)} e^{-t x^{-\gamma}} g_{\alpha}(t) dt \\ \label{30.10.13-3b}
& = \gamma x^{-(1+\gamma)} \left[-\frac{d}{du} \int_{0}^{\infty} e^{- t u} g_{\alpha}(t) dt\right]_{u=x^{-\gamma}} \\ \label{30.10.13-3c}
& = \gamma x^{-(1+\gamma)} \left[-\frac{d}{du} e^{-u^{\alpha}}\right]_{u=x^{-\gamma}} \\ \label{30.10.13-3d}
& = \gamma x^{-(1+\gamma)} \left[\alpha u^{\alpha-1} e^{-u^{\alpha}}\right]_{u=x^{-\gamma}} \\ \label{30.10.13-3e}
& = \gamma\alpha x^{-(1+\gamma\alpha)} e^{-x^{-\gamma\alpha}}, \quad \gamma > 0, 
\end{align}
so it is again given by the Fr\'{e}chet pdf. Therefore, the Laplace transform of $\bar{g}_{\alpha}(\gamma, x)$ is given by
\begin{equation}\label{30.10.13-4}
\mathcal{L}[\bar{g}_{\alpha}(\gamma, x); p] = \mathcal{L}[Fr(\gamma\alpha, x); p],
\end{equation}
and according to Eq. \eqref{eq8} can be exactly computed for rational $\gamma\alpha$. 

The Fr\'{e}chet transform of another Fr\'{e}chet distribution is equal to:
\begin{align}\label{30.10.13-5}
\bar{Fr}_{\gamma}(\beta, x) & = \int_{0}^{\infty} \sigma_{\gamma}(x, t) Fr(\beta, t) dt = \int_{0}^{\infty} \gamma t x^{-(1+\gamma)} e^{-tx^{-\gamma}} Fr(\beta, t) dt \nonumber\\
& = \gamma x^{-(1+\gamma)} \left\{-\frac{d}{du} \mathcal{L}[Fr(\beta, t); u]\right\}_{u=x^{-\gamma}},
\end{align}
and again can be calculated, for arbitrary positive $\gamma$ and rational $\beta$, from derivatives of Eq. \eqref{eq8}.

\begin{figure}[!h]
\begin{center}
\includegraphics[scale=0.56]{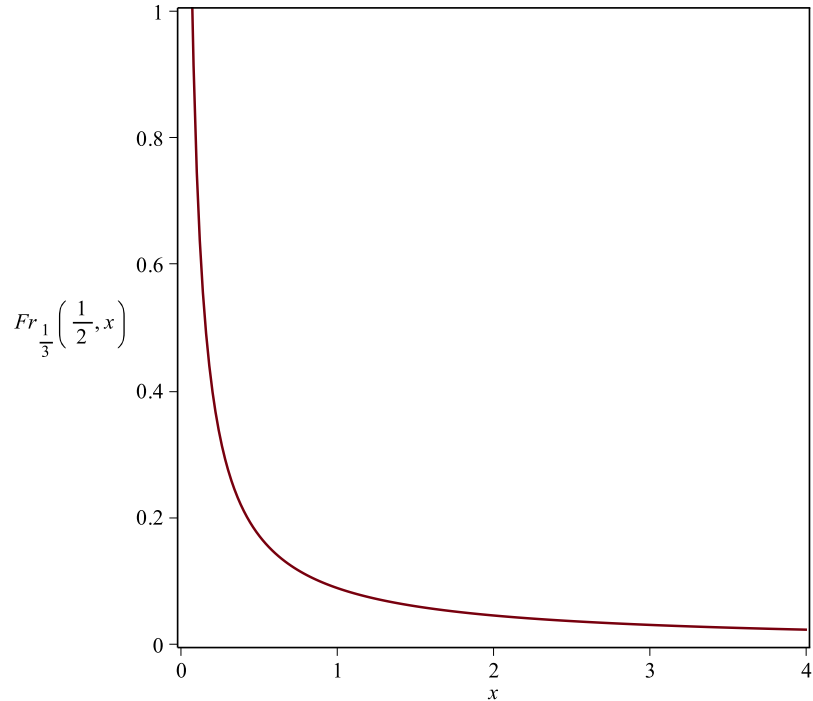}
\caption{\label{fig3}Plot of Eq. \eqref{16.02.2014-1} for $\gamma = 1/3$.}
\end{center}
\end{figure}

Since the derivatives of Meijer G functions are again expressible by Meijer G functions (see formula 8.2.2.32 of [\onlinecite{APPrudnikov-v3}]), the transform of type of Eq. \eqref{30.10.13-5} can be exactly evaluated. If we choose for illustration $\beta=1/2$, then for arbitrary $\gamma > 0$
\begin{align}\label{16.02.2014-1}
\begin{split}
\bar{Fr}_{\gamma}(1/2, x) = \int_{0}^{\infty} \sigma_{\gamma}(x, t) Fr(1/2, t) dt &= \ulamek{\gamma}{4\sqrt{\pi}} x^{-(1+\gamma)} MeijerG([[\,\,\,], [\,\,\,]], [[-\ulamek{1}{2}, 0, 0], [\,\,\,]], \ulamek{x^{-\gamma}}{4}) \\
& = \ulamek{\gamma}{4\sqrt{\pi}} x^{-(1+\gamma)} G^{3, 0}_{0, 3} \left(\ulamek{x^{-\gamma}}{4}\Big\vert{- \atop -\ulamek{1}{2}, 0, 0}\right).
\end{split}
\end{align}
The graphical representation of Eq. \eqref{16.02.2014-1} is given on Fig. \ref{fig3} for $\gamma = 1/3$.
We provide below for completeness the demonstration of Eq. \eqref{16.02.2014-1}:
\begin{align}\label{27Feb14-1}
\bar{Fr}(1, 2, x) &= \gamma x^{-(1+\gamma)} \left[-\frac{d}{du} \pi^{-1/2} G^{3, 0}_{0, 3}\left(\frac{u}{4}\Big\vert {- \atop 0, \ulamek{1}{2}, 1}\right)\right]_{u=x^{-\gamma}} \\ \label{27Feb14-2}
&= -\frac{\gamma}{4\sqrt{\pi}} x^{-(1+\gamma)} \left[\frac{d}{dv} G^{3, 0}_{0, 3}\left(v \Big\vert {- \atop 0, \ulamek{1}{2}, 1}\right)\right]_{v=x^{-\gamma}/4}\\ \label{27Feb14-3}
& = \frac{\gamma}{4\sqrt{\pi}} x^{-(1+\gamma)} \left[-v^{-1} G^{3, 1}_{1, 4}\left(v \Big\vert {0 \atop 0, \ulamek{1}{2}, 1, 1}\right)\right]_{v=x^{-\gamma}/4} \\ \label{27Feb14-4}
& = \frac{\gamma}{4\sqrt{\pi}} x^{-(1+\gamma)} \left[v^{-1} G^{4, 0}_{1, 4}\left(v \Big\vert {0 \atop 0, \ulamek{1}{2}, 1, 1}\right)\right]_{v=x^{-\gamma}/4} \\ \label{27Feb14-5}
&= \frac{\gamma}{4\sqrt{\pi}} x^{-(1+\gamma)} \left[v^{-1} G^{3, 0}_{0, 3}\left(v \Big\vert {- \atop \ulamek{1}{2}, 1, 1}\right)\right]_{v=x^{-\gamma}/4} \\ \label{27Feb14-6}
& = \frac{\gamma}{4\sqrt{\pi}} x^{-(1+\gamma)} \left[G^{3, 0}_{0, 3}\left(v \Big\vert {- \atop -\ulamek{1}{2}, 0, 0}\right)\right]_{v=x^{-\gamma}/4} \\ \label{27Feb14-7}
& = \frac{\gamma}{4\sqrt{\pi}} x^{-(1+\gamma)} G^{3, 0}_{0, 3}\left(\frac{x^{-\gamma}}{4} \Big\vert {- \atop -\ulamek{1}{2}, 0, 0}\right).
\end{align}
Eq. \eqref{27Feb14-3} has been obtained using the differentiation formula (8.2.2.32) of [\onlinecite{APPrudnikov-v3}], Eq. \eqref{27Feb14-4} results from formula (8.2.2.16) of [\onlinecite{APPrudnikov-v3}] and Eq. \eqref{27Feb14-6} is based on the translation formula (8.2.2.15) of [\onlinecite{APPrudnikov-v3}].

\section{Summary and Conclusions}

We are in a position now to confront several characteristics of Fr\'{e}chet and L\'{e}vy stable laws. The functional form of the Fr\'{e}chet pdf, see Eq. \eqref{eq1}, is an elementary function displaying the essential singularity at $x\to 0$, unimodality and heavy-tailed decay at $x\to\infty$ for all the values of $\gamma > 0$. The one-dimensional L\'{e}vy laws $g_{\alpha}(x)$ can be represented by elementary function exclusively for $\alpha = 1/2$. Their explicit forms can be expressed by standard special functions only for $\alpha = 1/3, 2/3$; otherwise for other values of rational $\alpha$ they can be represented by finite sum of generalized hypergeometric functions [\onlinecite{KAPenson10}]. For arbitrary real $0 < \alpha < 1$ the explicit form of $g_{\alpha}(x)$ is unknown. The Laplace transform of Fr\'{e}chet distribution, as shown in the present work, can be expressed by higher transcendental function, the Meijer G function, but exclusively for rational values of the shape parameter $\gamma$. The Laplace transform for general real $\gamma > 0$ is unknown. In contrast, the Laplace transform of L\'{e}vy stable law $g_{\alpha}(x)$ is given for any $0 < \alpha < 1$ by the stretched exponential $\exp(-p^{\alpha})$. 

\begin{figure}[!h]
\begin{center}
\includegraphics[scale=0.56]{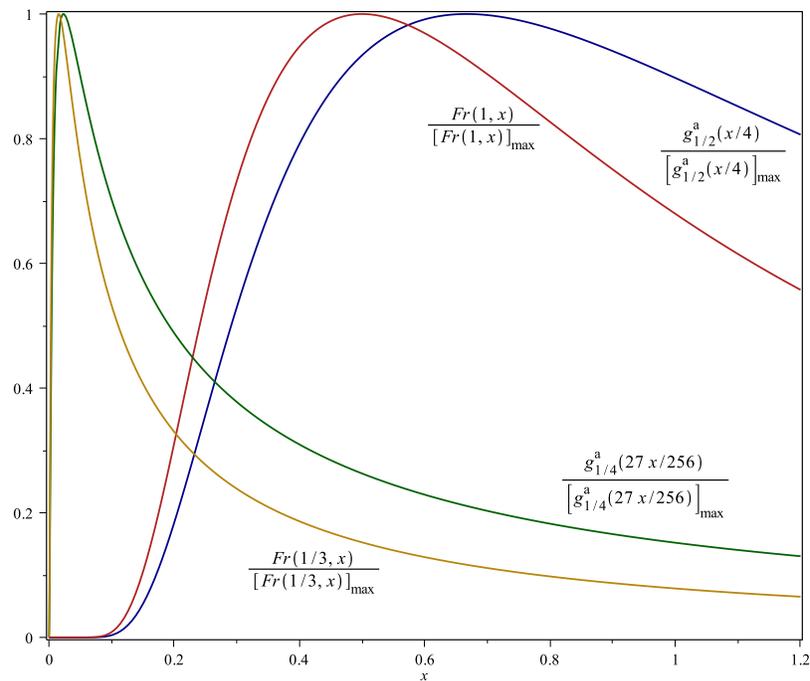}
\caption{\label{fig4} (Color online) Comparison of the reduced $g^{\rm a}_{\alpha}(t)$ with the reduced $Fr(\gamma, x)$ for $\alpha = 1/2$, $t=x/4$ (blue line), and $\gamma =1$ (red line), and for $\alpha = 1/4$, $t=27x/256$ (green line), and $\gamma=1/3$ (gold line).}
\end{center}
\end{figure}

This comparison can be carried even further by exploiting the asymptotic forms of $g_{\alpha}(x)$, $x\to 0$, obtained in [\onlinecite{Mikusinski}]. In fact the asymptotic formula
\begin{equation}\label{4Mar14-1}
g^{\rm a}_{\alpha}(t) = \frac{1}{\sqrt{2\pi}} \alpha^{\frac{1}{2-2\alpha}} t^{-\frac{2-\alpha}{2-2\alpha}} \exp\left[-(1-\alpha) \alpha^{\frac{\alpha}{1-\alpha}} t^{-\frac{\alpha}{1-\alpha}}\right],
\end{equation}
after introducing there $\alpha = \gamma/(1+\gamma)$ and $t = \gamma x/(1+\gamma)^{1+1/\gamma}$ furnishes
\begin{equation}\label{4Mar14-2}
g^{\rm a}_{\gamma/(1+\gamma)}\left(\frac{\gamma x}{(1+\gamma)^{1+1/\gamma}}\right) = \frac{(1+\gamma)^{1/\gamma}}{\sqrt{2\pi}} \left(\frac{1+\gamma}{\gamma}\right)^{3/2} x^{\gamma/2} Fr(\gamma, x).
\end{equation}
Therefore, for given $\gamma$, and for small $x$ it is convenient to compare $Fr(\gamma, x)$ with the right hand side of Eq. \eqref{4Mar14-2}. Spurious effects can be eliminated by considering reduced distributions, obtained by dividing the original quantities by their corresponding maximal values. This is done in Fig. \ref{fig4}. As expected from Eq. \eqref{4Mar14-2} the two distributions approach each other even more, as $\gamma$ decreases.

Both of these distributions have known expressions of their $\mu$-th power moment: for the L\'{e}vy laws the only non-diverging moments are equal to 
\begin{equation}\label{28Feb14-1}
\int_{0}^{\infty} x^{\mu} g_{\alpha}(x) dx = \frac{\Gamma(1-\ulamek{\mu}{\alpha})}{\Gamma(1-\mu)}, \quad -\infty < \mu < \alpha.
\end{equation}
The corresponding expression for the Fr\'{e}chet law is
\begin{equation}\label{28Feb14-2}
\int_{0}^{\infty} x^{\mu} Fr(\gamma, x) dx = \Gamma(1-\ulamek{\mu}{\gamma}), \quad -\infty < \mu < \gamma.
\end{equation}
Thus, the divergence of sufficiently high moments is a common feature of both of these distributions. 

For a deeper probabilistic analysis of the comparison between the L\'{e}vy and Fr\'{e}chet laws the reader is referred to [\onlinecite{TSimon14}].

Last but not least we would like to underline the utility of Meijer G functions which are the most important calculational tool in the present work. For reasons explored above our key results can not be represented by any other better known functions. In our opinion the Meijer G functions must be considered as known special functions whose diverse applications are becoming truly universal: for a choice of recent applications see [\onlinecite{DBermudez14, PJForrester14, WMlotkowski13, GAkemann13}] and the references therein. In a recent exposition [\onlinecite{RBeals13}] their use in strongly advocated as being ``both natural and attractive''. Needless to say, we unreservedly share this conviction.   

\section*{Acknowledgment}

The authors acknowledge support from the PHC Polonium, Campus France, project no. 28837QA.


\begin{thebibliography}{99}
\bibitem{CAS} We have made extensive use of Mathematica$^{\text{\textregistered}}$ and Maple$^{\text{\textregistered}}$ in this work.
\end{thebibliography}
\end{document}